\documentclass[12pt]{article}
\usepackage{amscd,amsmath,amssymb,amsfonts,times,ascmac}
\usepackage{mathrsfs}
\usepackage[dvips]{color} 
\usepackage[all]{xy}
\usepackage{enumerate}
 \usepackage{fancybox} 

\usepackage{amsmath}
\usepackage{amssymb}

\numberwithin{equation}{section}
    \newtheorem{thm}{Theorem}[section]
    \newtheorem{lem}[thm]{Lemma}
    
    \newtheorem{prop}[thm]{Proposition}
    \newtheorem{cor}[thm]{Corollary}

    \newtheorem{defn}[thm]{Definition}
    
    \newtheorem{rem}[thm]{Remark}
\DeclareMathAlphabet{\mathpzc}{OT1}{pzc}{m}{it}
\newcommand{\qed}
{\mbox{}\nolinebreak$\square$\medbreak\par}
\newenvironment{pf}{\par\smallskip\noindent\emph{Proof.}}{\hfill\qed\par\smallskip}
\newenvironment{pf*}[1]{\par\smallskip\noindent\emph{#1.}}{\hfill\qed\par\smallskip}
\newcommand{\bysame}{\hskip.3em \leavevmode\rule[.5ex]{2.5em}{.3pt}\hskip0.5em}
\begin{document}
\title{Zeta functions of certain K3 families : application of the formula of Clausen}
\author{M. Asakura\footnote{
Hokkaido University, 
Sapporo 060-0810, JAPAN. 
\texttt{asakura@math.sci.hokudai.ac.jp}
}
}
\date\empty
\maketitle

\def\can{\mathrm{can}}
\def\cano{\mathrm{canonical}}
\def\ch{{\mathrm{ch}}}
\def\Coker{\mathrm{Coker}}
\def\crys{\mathrm{crys}}
\def\dlog{{\mathrm{dlog}}}
\def\dR{{\mathrm{d\hspace{-0.2pt}R}}}            
\def\et{{\mathrm{\acute{e}t}}}  
\def\Frac{{\mathrm{Frac}}}
\def\phami{\phantom{-}}
\def\id{{\mathrm{id}}}              
\def\Image{{\mathrm{Im}}}        
\def\Hom{{\mathrm{Hom}}}  
\def\Ext{{\mathrm{Ext}}}
\def\MHS{{\mathrm{MHS}}}  
  
\def\can{\mathrm{can}}
\def\Tr{\mathrm{Tr}}
\def\ch{{\mathrm{ch}}}
\def\Coker{\mathrm{Coker}}
\def\crys{\mathrm{crys}}
\def\dlog{d{\mathrm{log}}}
\def\dR{{\mathrm{d\hspace{-0.2pt}R}}}            
\def\et{{\mathrm{\acute{e}t}}}  
\def\Frac{\operatorname{Frac}}
\def\phami{\phantom{-}}
\def\id{{\mathrm{id}}}              
\def\Image{{\mathrm{Im}}}        
\def\Hom{\operatorname{Hom}}  
\def\Ext{{\mathrm{Ext}}}
\def\MHS{{\mathrm{MHS}}}  
  
\def\ker{\operatorname{Ker}}          
\def\mf{{\text{mapping fiber of}}}
\def\Pic{{\mathrm{Pic}}}
\def\CH{{\mathrm{CH}}}
\def\NS{{\mathrm{NS}}}
\def\NF{{\mathrm{NF}}}
\def\End{\operatorname{End}}
\def\pr{{\mathrm{pr}}}
\def\red{{\mathrm{red}}}
\def\Proj{\operatorname{Proj}}
\def\ord{\operatorname{ord}}
\def\rig{{\mathrm{rig}}}
\def\reg{{\mathrm{reg}}}          %
\def\res{{\mathrm{res}}}          %
\def\Res{\operatorname{Res}}
\def\Spec{\operatorname{Spec}}     
\def\syn{{\mathrm{syn}}}
\def\Gr{{\mathrm{Gr}}}
\def\ln{{\operatorname{ln}}}

\def\bA{{\mathbb A}}
\def\bC{{\mathbb C}}
\def\C{{\mathbb C}}
\def\G{{\mathbb G}}
\def\bE{{\mathbb E}}
\def\bF{{\mathbb F}}
\def\F{{\mathbb F}}
\def\bG{{\mathbb G}}
\def\bH{{\mathbb H}}
\def\bJ{{\mathbb J}}
\def\bL{{\mathbb L}}
\def\cL{{\mathscr L}}
\def\bN{{\mathbb N}}
\def\bP{{\mathbb P}}
\def\P{{\mathbb P}}
\def\bQ{{\mathbb Q}}
\def\Q{{\mathbb Q}}
\def\bR{{\mathbb R}}
\def\R{{\mathbb R}}
\def\bZ{{\mathbb Z}}
\def\Z{{\mathbb Z}}
\def\cH{{\mathscr H}}
\def\cD{{\mathscr D}}
\def\cE{{\mathscr E}}
\def\cF{{\mathscr F}}
\def\cU{{\mathscr U}}
\def\O{{\mathscr O}}
\def\cR{{\mathscr R}}
\def\cS{{\mathscr S}}
\def\cC{{\mathscr C}}
\def\cX{{\mathscr X}}
\def\cY{{\mathscr Y}}
%
\def\ve{\varepsilon}
\def\vG{\varGamma}
\def\vg{\varGamma}
%
%
%
%
\def\lra{\longrightarrow}
\def\lla{\longleftarrow}
\def\Lra{\Longrightarrow}
\def\hra{\hookrightarrow}
\def\lmt{\longmapsto}
\def\ot{\otimes}
\def\op{\oplus}
\def\l{\lambda}
\def\Isoc{{\mathrm{Isoc}}}
\def\Fil{{\mathrm{Fil}}}

\def\MHS{{\mathrm{MHS}}}
\def\Dw{{\mathrm{Dw}}}
\def\VMHS{{\mathrm{VMHS}}}
\def\sm{{\mathrm{sm}}}
\def\tr{{\mathrm{tr}}}
\def\fib{{\mathrm{fib}}}
\def\FIsoc{{F\text{-Isoc}}}
\def\FMIC{{F\text{-MIC}}}
\def\Log{{\mathscr{L}{og}}}
\def\FilFMIC{{\mathrm{Fil}\text{-}F\text{-}\mathrm{MIC}}}

\def\wt#1{\widetilde{#1}}
\def\wh#1{\widehat{#1}}
\def\spt{\sptilde}
\def\ol#1{\overline{#1}}
\def\ul#1{\underline{#1}}
\def\us#1#2{\underset{#1}{#2}}
\def\os#1#2{\overset{#1}{#2}}

\def\eq{\quad\Longleftrightarrow\quad}
\def\Ross{{\xi_{\mathrm{Ross}}}}
\def\HG{{\mathrm{HG}}}
\def\Gal{{\mathrm{Gal}}}
\def\rank{{\mathrm{rank}}}
\def\Sym{{\mathrm{Sym}}}
\def\Frob{{\mathrm{Frob}}}
\def\a{\alpha}
\begin{abstract}
Based on the theory of rigid cohomology, we provide an explicit formula of zeta functions of
certain K3 families, which we call the hypergeometric type.
The central point of our argument is the comparison between the 2nd rigid cohomology of a K3 and
the symmetric product of an elliptic curve, that is brought from the classical
formula of Clausen.
\end{abstract}

\section{Introduction}\label{intro-sect}
A projective smooth surface $X$ is called a K3 surface if
\[
H^1(X,\O_X)=0,\quad K_X\cong \O_X.
\]
The subject of this paper is the $l$-adic Galois representation
\[
G_\Q=\Gal(\ol\Q/\Q)\lra \mathrm{Aut}(H^2_\et(\ol X,\Q_l)),\quad \ol X:=X\times_\Q\ol\Q
\]
for a K3 surface $X$ over $\Q$ such that the Picard number $\rho(\ol X)
:=\rank\NS(\ol X)$ is $\geq 19$.
Thanks to the theorem of Morrison \cite{morrison},
$\ol X$ is isogenous
to a Kummer K3 surface $\mathrm{Km}(\ol E\times \ol E)$ of an elliptic curve $E$ 
(not necessarily defined over 
$\Q$),
which is often referred to as the Shioda-Inose structure.
Then one finds that the Galois representation $H^2_\et(\ol X,\Q_l)$ is
potentially isomorphic to the symmetric product $\Sym^2H^1_\et(\ol E,\Q_l)$
up to a simple factor.
Morrison's theorem asserts only the existence of the Kummer K3,
the problem on finding an explicit $E$ is nontrivial.
Besides, the isogeny is not necessarily defined over $\Q$ (even when so is $E$),
and then it is another new task to explore the $G_\Q$-representation.

\medskip

In this paper, we study
the characteristic polynomial
\begin{equation}\label{intro-eq1}
\det(1-\phi_p^{-1}T\mid H^2_\et(\ol X,\Q_l))
\end{equation}
of the $p$-th Frobenius $\phi_p\in G_\Q$
for a K3 surface $X$ whose ``period" is
the hypergeometric series
\[
F_{\ul \a}(t)={}_3F_2\left({\a_0,\a_1,\a_2\atop 1,1};t\right)=\sum_{n=0}^\infty
\frac{(\alpha_0)_n}{n!}\frac{(\alpha_1)_n}{n!}\frac{(\alpha_2)_n}{n!}t^n,
\quad (\alpha)_n:=\alpha(\alpha+1)\cdots(\alpha+n-1)
\]
where 
$\ul\a=(\a_0,\a_1,\a_2)$ is either of the following,
\begin{equation}\label{alpha}
\left(\frac12,\frac12,\frac12\right),\left(\frac13,\frac23,\frac12\right),
\left(\frac14,\frac34,\frac12\right),
\left(\frac16,\frac56,\frac12\right).
\end{equation}
In precise, we consider a projective smooth family
\[
f:\cX\lra T=\Spec \Q[t,(t-t^2)^{-1}]
\]
of K3 surfaces such that the generic fiber 
$\ol X_t=\cX\times_T\ol{\Q(t)}$ satisfies $\rho(\ol X_t)=19$.
Put
\[
V_\dR(\cX/T)=\mathrm{Coim}[H^2_\dR(\cX/T)\to H^2_\dR(\ol X_t)/\NS(\ol X_t)\ot\ol{\Q(t)}]
\]
a free $\O(T)$-module of rank $3$.
Let $\cD=\Q[t,(t-t^2)^{-1},\frac{d}{dt}]$ be the Wyle algebra of $T$, and let
\[
P_{\ul\a}=D^3-t(D+\a_0)(D+\a_1)(D+\a_2),\quad D:=t\frac{d}{dt}
\]
be the hypergeometric differential operator, which annihilates $F_{\ul\a}(t)$.
Then we call $f$ {\it of hypergeometric type} $F_{\ul\a}(t)$ if there is an isomorphism
\[
V_\dR(\cX/T)\cong \cD/\cD P_{\ul\a}
\]
of $\cD$-modules (Definition \ref{K3.HG}).
Here are examples of K3 families of hypergeometric type.
\begin{itemize}
\item[(i)]
The Dwork family (cf. \cite{Katz})
\begin{equation}\label{intro-ex1}
tx_0^4+x_1^4+x_2^4+x_3^4-4x_0x_1x_2x_3=0
\end{equation}
of quartic surfaces is of hypergeometric type $F_{\frac14,\frac34,\frac12}(t)$.
\item[(ii)]
The K3 family (cf.  \cite{AOP})
\begin{equation}\label{intro-ex2}
z^2=xy(1+x)(1+y)(x-t y)
\end{equation}
is of hypergeometric type $F_{\frac12,\frac12,\frac12}(t)$.
\item[(iii)]
The K3 family (cf. \cite[\S 6.3]{As-Ross1})
\begin{equation}\label{intro-ex3}
(1-x^2)(1-y^2)(1-z^2)=t
\end{equation}
is of hypergeometric type $F_{\frac12,\frac12,\frac12}(t)$.
This is isogenous
to the family \eqref{intro-ex2} over $\Q$ (\cite[Lemma 6.3]{As-Ross1}).
\item[(iv)]
Let $n=3,4,6$, and $\cE_n^\pm\to T$ the elliptic K3 surface
constructed in
\cite[6.4]{As-Ross2}.
Then this is of hypergeometric type $F_{\frac1n,\frac{n-1}{n},\frac12}(t)$.
\end{itemize}

The purpose of this paper is that for each $\ul\a$ in \eqref{alpha},
we describe the characteristic polynomial \eqref{intro-eq1} 
by a specific elliptic curve 
\begin{equation}\label{Ea-defn}
E_{\ul\a,s}=
\begin{cases}
y^2=x(x-1)(x-s)&\ul\a=(\frac12,\frac12,\frac12)\\
y^2=x^3+(3x+4-4s)^2&\ul\a=(\frac13,\frac23,\frac12)\\
y^2=x(x^2-2x+1-s)&\ul\a=(\frac14,\frac34,\frac12)\\
y^2=4x^3-3x+1-2s&\ul\a=(\frac16,\frac56,\frac12).
\end{cases}
\end{equation}

\begin{thm}[Theorem \ref{thm-1}]\label{intro-thm-1}
Let $f:\cX\to T$ be a K3 family of hypergeometric type $F_{\ul\a}(t)$.
Let $p>3$ be a prime at which  there is an integral regular flat model
\[
f_{\Z_{(p)}}:\cX_{\Z_{(p)}}\lra T_{\Z_{(p)}}
\]
over the ring $\Z_{(p)}\subset \Q$ such that $f_{\Z_{(p)}}$ is smooth projective. 
Let $a\in \Z_p$ such that $a(1-a)\not\equiv0$ mod $p$, and $X_a$ the fiber
at $t=a$.
Let 
\[
V_\et(\ol X_a)_{\Q_l}:=\mathrm{Coim}[H^2_\et(\ol X_a,\Q_l)
\to H^2_\et(\ol X_t,\Q_l)/\NS(\ol X_t)\ot\Q_l]\cong \Q_l^3.
\]
Put $b=\frac12(1-\sqrt{1-a})$, and let $E_{\ul\a,b}$ be the elliptic curve \eqref{Ea-defn}.
Let $1-a_{p^2}(E_{\ul\a,b})T+p^2T^2$ be the characteristic polynomial of 
the $p^2$-th Frobenius
$\phi_{p^2}\in G_\Q$, namely $a_{p^2}(E_{\ul\a,b})\in \Z$ satisfies
\[
1-a_{p^2}(E_{\ul\a,b})+p^2=\sharp E_{\ul\a,b}(\F_{p^2}).
\]
Let
\[
d_{\ul\a}=\begin{cases}
-1&\ul\a=(\frac12,\frac12,\frac12),(\frac16,\frac56,\frac12)\\
-2&\ul\a=(\frac14,\frac34,\frac12)\\
-3&\ul\a=(\frac13,\frac23,\frac12)
\end{cases}
\]
and put
\[
A_{a,p}=\begin{cases}
a_{p^2}(E_{\ul\a,b})&\sqrt{1-a}\in \Z_p\\
(\frac{d_{\ul\a}}{p})a_{p^2}(E_{\ul\a,b})
&\sqrt{1-a}\not\in \Z_p,\,E_{\ul\a,b}\text{: ordinary at $p$}\\
2p&\sqrt{1-a}\not\in \Z_p,\,E_{\ul\a,b}\text{: supersingular at $p$}.
\end{cases}
\]
where $(\frac{\ast}{p})$ denotes the Legendre symbol.
Then
\[
\det(1-\phi_p^{-1}T\mid V_\et(\ol X_a)_{\Q_l})
=\left(1-\left(\frac{1-a}{p}\right)
\chi_{\cX/T}(\phi_p)pT\right)(1-\chi_{\cX/T}(\phi_p)A_{a,p}T+p^2T^2)
\]
where 
$\chi_{\cX/T}:G_\Q\to\{\pm1\}$ is the character \eqref{chiX}
defined in \S \ref{K3HG-sect}.
\end{thm}
Notice that $A_{a,p}$ does not depend on the choice of $b$ 
as $E_{\ul\a,b}$ and $E_{\ul\a,1-b}$ are isogenous over $\F_{p}(b,\sqrt{d_{\ul\a}})$
(see \eqref{isogeny1}, \ldots,\eqref{isogeny4} below).
For the families (i), \ldots, (iv), the character $\chi_{\cX/T}$ is trivial 
(Remark \ref{chi-rem}).

\medskip

As a byproduct of the proof of Theorem \ref{intro-thm-1},
we have the following description of the Galois representation over $\Q(\sqrt{1-a})$. 
\begin{cor}[Theorem \ref{thm-2}]\label{intro-cor}
Let $F$ be a number field and 
let $a\in F\setminus\{0,1\}$ be arbitrary.
Then
there is an isomorphism
\begin{equation}\label{intro-eq2}
V_\et(\ol X_a)_{\Q_l}\cong \Sym^2H^1_\et(\ol E_{\ul\a,b},\Q_l)\ot\chi_{\cX/T}
\end{equation}
of $G_{F(\sqrt{1-a})}$-representations. 
\end{cor}
If
$\ul\a=(\frac12,\frac12,\frac12)$ or $(\frac16,\frac56,\frac12)$,
we have an alternative description of the characteristic polynomial
by another elliptic curve,
\begin{equation}\label{Ca-defn}
C_{\ul\a,t}=
\begin{cases}
y^2=x^3-2x^2+\frac{t}{t-1}x
&\ul\a=(\frac12,\frac12,\frac12)\\
y^2=4x^3-3(1-t)x+(1-t)^2
&\ul\a=(\frac16,\frac56,\frac12).
\end{cases}
\end{equation}

\begin{thm}[Theorem \ref{thm-3}]\label{intro-thm-2}
Let $\ul\a$ be either of $(\frac12,\frac12,\frac12)$ or $(\frac16,\frac56,\frac12)$.
Then
\begin{align*}
&\det(1-\phi_p^{-1}T\mid V_\et(\ol X_a)_{\Q_l})\\
=&\left(1-\left(\frac{1-a}{p}\right)\chi_{\cX/T}(\phi_p)pT\right)
\left(1-\left(\frac{1-a}{p}\right)\chi_{\cX/T}(\phi_p)a_{p^2}(C_{\ul\a,a})T+p^2T^2\right).
\end{align*}
Hence, for arbitrary $a\in \Q\setminus\{0,1\}$, there is an isomorphism
\begin{equation}\label{intro-eq3}
V_\et(\ol X_a)_{\Q_l}\cong \Sym^2H^1_\et(\ol C_{\ul\a,b},\Q_l)\ot
\chi_{\cX/T}\ot\chi_{1-a}
\end{equation}
of $G_\Q$-representations where
$\chi_{1-a}$ denotes the Kronecker character for $\Q(\sqrt{1-a})$.
\end{thm}

\medskip

There are lots of works concerning zeta functions of K3 surfaces or
Calabi-Yau manifolds with hypergeometric functions, 
the author does not catch up all of them though. 
For instance, there are a number of papers describing the zeta functions in terms of finite hypergeometric
series, \cite{Go1}, \cite{Go2}, \cite{Koblitz}, \cite{Mc}, \cite{Miyatani}, \cite{otsubo} etc.
On the other hand, the author finds only a few papers which
exhibit a characteristic polynomial of $\phi_p$ (not $\phi_{p^m}$ for a particular $m$)
for all but finitely many $p$.
Concerning the Dwork family (i), the Shioda-Inose structure
is provided by Elkies-Sch\"utt \cite{ES}, which imposes the Galois representation
$V_\et(\ol X_a)_{\Q_l}$ potentially.
Corollary \ref{intro-cor} for $X$ the Dwork family
can be derived from Naskrecki \cite[Cor.6.7]{nas} or
Otsubo \cite[Thm.7.4]{otsubo} where they discuss the family in a context of finite hypergeometric series.
However these results are not enough to determine
the $G_\Q$-representation. As long as the author sees, Theorem \ref{intro-thm-1} 
is a new formula for the Dwork family.
Concerning the K3 family (ii),
the isomorphism \eqref{intro-eq3} is proved in \cite{AOP},
and the Shioda-Inose structure (defined over $\Q$)
is exhibited in \cite{GT}.
For this family,
Theorem \ref{intro-thm-2} 
is nothing new, while our proof is entirely different from theirs.

\medskip

For the proof of Theorems \ref{intro-thm-1} and \ref{intro-thm-2},
we follow the argument of Dwork \cite{Dwork-unique}.
The key tool is the rigid cohomology (the Monsky-Washnitzer cohomology) (cf. \cite{LS}).
The characteristic polynomial of Frobenius
can be obtained from the Frobenius structure on the rigid cohomology
\[
H_\rig^2(\cX_{\F_p}/T_{\F_p})\]
where $\cX_{\F_p}:=\cX_{\Z_{(p)}}\times_{\Z_{(p)}}\F_p$ etc.
We then compare it with the symmetric product
\[
\Sym^2H^1_\rig(\cE_{\ul\a,\F_p}/S_{\F_p})
\]
in a direct way, where 
$S=\Spec\Q[t,\sqrt{1-t},(t-t^2)^{-1}]$ 
and $\cE_{\ul\a}\to S$ is the family of elliptic curves $E_{\ul\a,s}$
with $s=\frac12(1-\sqrt{1-t})$.
See Theorem \ref{ell-F-lem2} for the detail. 
The comparison is brought from the classical formula of Clausen (\cite[16.12.2]{NIST})
\begin{equation}\label{Clausen}
{}_3F_2\left({2a,2b,a+b\atop a+b+\frac12,2a+2b};t\right)=
{}_2F_1\left({a,b\atop a+b+\frac12};t\right)^2.
\end{equation}
The idea is sketched in \cite[p.92--93]{Dwork-unique} for the Dwork family, 
while we take more thorough discussion in this paper.
Besides, to work out on the ``sign" such as $(\frac{d_{\ul\a}}{p})$,
we need additional argument that is not suggested in loc.cit.
As a final comment, Otsubo's approach is comparable with ours. 
He obtains an analogue of Clausen's formula
in a context of finite hypergeometric series (\cite[Thm. 6.5]{otsubo}), and proves a similar (but weaker) result to Theorem \ref{intro-thm-1}
for the Dwork family (\cite[Thm. 7.4]{otsubo}).

\medskip

\noindent\textbf{Acknowledgement}.
The author is grateful to Noriyuki Otsubo for the stimulating discussion on 
the Dwork family 
and for 
 encouraging him to write this paper.

\section{Characteristic polynomial for Hypergeometric differential equations }\label{Fisoc-sect}
Let $p$ be a prime number. Let $W=W(\ol\F_p)$ be the Witt ring of the algebraic closure
$\ol\F_p$ of $\F_p$.
Let $K=\Frac(W)$ be the fractional field.
\subsection{$F$-isocrystals of hypergeometric differential equations}\label{F-H-sect}
For $\a\in \Z_p$, we denote by $\a'$ the Dwork prime which is defined to be
$(\a+l)/p$ where $l\in\{0,1,\ldots,p-1\}$ such that $a+l\equiv 0$ mod $p$.
We define $\a^{(i)}:=(\a^{(i-1)})^\prime$ with $\a^{(0)}=\a$.
For $\ul\a=(\a_0,\a_1,\ldots,\a_n)\in\Z_p^{n+1}$, we denote
$\ul \a^{(i)}=(\a^{(i)}_0,\a^{(i)}_1,\ldots,\a^{(i)}_2)$.
We write
\[
F_{\ul\a}(t)
={}_{n+1}F_n\left({\a_0,\ldots,\a_n\atop 1,\ldots,1};t\right)
=\sum_{i=0}^\infty\frac{(\a_0)_i}{i!}\cdots\frac{(\a_n)_i}{i!}t^i\in \Z_p[[t]]
\]
the hypergeometric series where $(\a)_i=\a(\a+1)\cdots(\a+i-1)$ is 
the Pochhammer symbol, and
\begin{equation}\label{PHG}
P_{\ul \a}
=D^{n+1}-t(D+\a_0)(D+\a_1)\cdots(D+\a_n),\quad 
D:=t\frac{d}{dt}
\end{equation}
the hypergeometric differential operator whose solution is $F_{\ul\a}(t)$.

\medskip

Let $n=2$ and $\ul \a=(\a_0,\a_1,\a_2)\in \Z_p^3$.
Let $\cD=K[t,(t-t^2)^{-1},\frac{d}{dt}]$ be the Wyle algebra, and let
\begin{equation}\label{VHG}
V_{\ul \a}:=\cD/\cD P_{\ul \a}=\cD\,\omega_{\ul \a}
\end{equation}
be a left $\cD$-module where the symbol $\omega_{\ul \a}$ denotes 
$1_\cD+\cD P_{\ul\a}$ 
a generator of $V_{\ul\a}$.
The module $V_{\ul \a}$ is a free $K[t,(t-t^2)^{-1}]$-module of rank $3$,
\[
V_{\ul \a}=\bigoplus_{i=0}^2K[t,(t-t^2)^{-1}]D^i\omega_{\ul \a}.
\]
We denote by $\nabla$ the connection on $V_{\ul\a}$ induced by $\frac{d}{dt}$.
Put
\[
y_2=(1-t)F_{\check{\ul \a}}(t),\quad y_1=-(\a_0+\a_1+\a_2)\frac{t}{1-t}y_2-Dy_2,\quad y_0=-(\a_0\a_1+\a_1\a_2+\a_0\a_2)\frac{t}{1-t}y_2-Dy_1
\]
and letting $\check{\ul \a}=(1-\a_0,1-\a_1,1-\a_2)$,
\begin{align}
\wh\omega_{\ul \a}&=\frac{1}{F_{\ul \a}(t)}\omega_{\ul \a}\label{F-H-eq1}\\
\wh\xi_{\ul \a}&=-(1-t)^{\a_0+\a_1+\a_2-1}\frac{1}{F_{\check{\ul \a}}(t)}
(DF_{\ul \a}(t)\cdot\omega_{\ul \a}-F_{\ul \a}(t)\cdot D\omega_{\ul \a})\\
\wh\eta_{\ul \a}&=y_0\omega_{\ul \a}+y_1D\omega_{\ul \a}+y_2D^2\omega_{\ul \a}.
\end{align}
Since
\begin{align*}
\omega_{\ul \a}&=F_{\ul \a}\wh\omega_{\ul \a}\\
D\omega_{\ul \a}&=
DF_{\ul \a}\cdot\wh\omega_{\ul \a}+
(1-t)^{1-\a_0-\a_1-\a_2}\frac{F_{\check{\ul \a}}}{F_{\ul \a}}\wh\xi_{\ul \a}\\
D^2\omega_{\ul \a}
&=
-y_2^{-1}(y_0F_{\ul \a}+y_1DF_{\ul \a})\wh\omega_{\ul \a}
-(1-t)^{-\a_0-\a_1-\a_2}\frac{y_1}{F_{\ul \a}}\wh\xi_{\ul \a}
+y_2^{-1}\wh\eta_{\ul \a}
\end{align*}
it turns out that $\{\wh\omega_{\ul \a},\wh\xi_{\ul \a},\wh\eta_{\ul \a}\}$ forms a
free basis of $W[[t]]\ot_{W[t]}V_{\ul \a}$.
Let 
\[
G_{\ul \a}(t):=(1-t)^{1-\a_0-\a_1-\a_2}F_{\ul \a}(t)^{-2}F_{\check{\ul \a}}(t),\quad
G_{\check{\ul \a}}(t):=(1-t)^{\a_0+\a_1+\a_2-2}F_{\check{\ul \a}}(t)^{-2}F_{\ul \a}(t),
\]
then it is straightforward to see
\begin{equation}\label{FrobX-eq1}
D\wh\omega_{\ul \a}
=G_{\ul \a}(t)\wh\xi_{\ul \a},\quad
D\wh\xi_{\ul \a}
=G_{\check{\ul \a}}(t)\wh\eta_{\ul \a},\quad D\wh\eta_{\ul \a}=0.
\end{equation}

\medskip

Let $\sigma$ be the $p$-th Frobenius on $W[[t]]$ given by $\sigma(t)=ct^p$
with $c\in 1+pW$.
Letting $A^\dag$ denote the weak completion of a ring $A$ (e.g. \cite[p.135]{LS}),
we put 
\begin{equation}\label{dagger-mod}
V_{\ul\a}^\dag:=K[t,(t-t^2)^{-1}]^\dag\ot_{K[t,(t-t^2)^{-1}]}V_{\ul\a}.
\end{equation}
\begin{lem}
Suppose $\a_i\in \Z_p\cap \Q$ and $0<\a_i<1$ for each $i$. Let $m\geq 1$ be the least integer such that
$\ul\a^{(m)}$ and $\ul\a$ agree up to permutation.
Then there is unique $\sigma$-linear $p$-th Frobenius endomorphism $\Phi$ on 
$\bigoplus_{i=0}^{m-1}V_{\ul\a^{(i)}}^\dag$ satisfying
$p\Phi D=D\Phi$,
$\Phi(V_{\ul\a^{(i)}}^\dag)\subset V_{\ul\a^{(i-1)}}^\dag$
and
\begin{equation}\label{FrobX-eq6}
\Phi(\wh\omega_{\ul\a^{(i)}})\equiv p^2\wh\omega_{\ul\a^{(i-1)}}\mod
\langle\wh\xi_{\ul\a^{(i-1)}},\wh\eta_{\ul\a^{(i-1)}}\rangle.
\end{equation}
\end{lem}
\begin{pf}
The uniqueness follows from \cite{Dwork-unique}.
The existence follows from the fact that there is a motive corresponding to $\op
V_{\ul\a^{(i)}}^\dag$.
Such a motive is given in \cite{As-Ross2}.
Let
\[
U=\Spec W[t,(t-t^2)^{-1}][x_0,x_1,x_2]/((1-x_0^{n_0})(1-x_1^{n_1})(1-x_2^{n_2})-t)
\]
be an affine scheme over $T=\Spec W[t,(t-t^2)^{-1}]$.
Then for $0<i_k<n_k$, 
\[
V_{\frac{i_0}{n_0},\frac{i_1}{n_1},\frac{i_2}{n_2}}\cong 
W_2H^2_\dR(U_K/T_K)(i_0,i_1,i_2)
\]
as $\cD$-module, and the Frobenius on 
\[
\op W_2H^2_\rig(U_{\ol\F_p}/T_{\ol\F_p})(i_0,i_1,i_2)
\]
satisfies the conditions in the lemma (\cite[Theorems 3.3, 4.5]{As-Ross2}).
\end{pf}
\begin{defn}\label{Vcrys}
We define 
 an $F$-isocrystal
\[
V^\crys(\ul \a,\sigma)=\left(\bigoplus_{i=0}^{m-1}V_{\ul\a^{(i)}},
\bigoplus_{i=0}^{m-1}V_{\ul\a^{(i)}}^\dag,
\nabla,\Phi,\sigma\right).
\]
Note $V^\crys(\ul \a,\sigma)=V^\crys(\ul \a^{(i)},\sigma)$.
Define $\Frob_{\ul\a}^{(\sigma)}$ to be the $3\times3$-matrix
defined by
\[
\begin{pmatrix}
\Phi(\omega_{\ul\a^{(1)}})&
\Phi(D\omega_{\ul\a^{(1)}})&
\Phi(D^2\omega_{\ul\a^{(1)}})
\end{pmatrix}
=\begin{pmatrix}
\omega_{\ul\a}&
D\omega_{\ul\a}&
D^2\omega_{\ul\a}
\end{pmatrix}\Frob_{\ul\a}^{(\sigma)}.
\]
\end{defn}
By a simple computation,
it is possible to describe the matrix $\Frob_{\ul\a}^{(\sigma)}$ explicitly except the
``constant terms". Let
\begin{align}\label{FrobX-eq3}
&\Phi(\wh\omega_{\ul\a^{(1)}})=p^2\wh\omega_{\ul\a}
-p^2E_1(t)\wh\xi_{\ul\a}+p^2E_2(t)\wh\eta_{\ul\a},\\
&\Phi(\wh\xi_{\ul\a^{(1)}})=p\epsilon\wh\xi_{\ul\a}
-pE_3(t)\wh\eta_{\ul\a}
\label{FrobX-eq4}\\
&\Phi(\wh\eta_{\ul\a^{(1)}})=\epsilon'\wh\eta_{\ul\a}
\label{FrobX-eq5}
\end{align}
with $\epsilon,\epsilon'\in K$ and $E_i(t)\in K[[t,t^{-1}]]$.
Apply $D$ on \eqref{FrobX-eq3}.
It follows from \eqref{FrobX-eq1} that we have
\[
pG_{\ul\a^{(1)}}(t^\sigma)\Phi(\wh\xi_{\ul\a})
=p^2(G_{\ul\a}(t)-tE'_1(t))\wh\xi_{\ul\a}\mod\langle\wh\eta_{\ul\a}\rangle.
\]
This implies $\epsilon =1$ by \eqref{FrobX-eq4}, and
\[
tE'_1(t)=
G_{\ul\a}(t)-G_{\ul\a^{(1)}}(t^\sigma)
\]
which characterizes the series expansion of $E_1(t)$ except the constant term.
Apply $D$ on \eqref{FrobX-eq4}.
We have
\[
pG_{\check{\ul \a}^{(1)}}(t^\sigma)\Phi(\wh\eta_{\ul\a})
=p(G_{\check{\ul \a}}(t)-tE'_3(t))\wh\eta_\cX,
\]
which implies $\epsilon'=1$ by \eqref{FrobX-eq5} and
\[
tE'_3(t)=
G_{\check{\ul \a}}(t)-G_{\check{\ul \a}^{(1)}}(t^\sigma).
\]
We again apply $D$ on \eqref{FrobX-eq3}.
We have
\[
tE'_2(t)=G_{\check{\ul \a}}(t)E_1(t)-G_{\ul\a^{(1)}}(t^\sigma)E_3(t).
\]
This implies $E_1(0)=E_3(0)$.
Summing up the above, we have 
\begin{prop}\label{Frob-K3}
We have $\epsilon=\epsilon'=1$
and $E_i(t)\in K[[t]]$ satisfy
\begin{align*}
tE'_1(t)&=
G_{\ul\a}(t)-G_{\ul\a^{(1)}}(t^\sigma)\\
tE'_3(t)&=
G_{\check{\ul \a}}(t)-G_{\check{\ul \a}^{(1)}}(t^\sigma)\\
tE'_2(t)&=G_{\check{\ul \a}}(t)E_1(t)-G_{\ul\a^{(1)}}(t^\sigma)E_3(t)
\end{align*}
and $E_1(0)=E_3(0)$.
\end{prop}
\subsection{Formula on Characteristic polynomial of $\Frob_{\ul\a}^{(\sigma)}$}
\label{Fisoc-formula-sect}
In what follows, let $\ul\a$ be either of the following,
\begin{equation}\label{F-for-eq1}
\left(\frac12,\frac12,\frac12\right),\left(\frac13,\frac23,\frac12\right),
\left(\frac14,\frac34,\frac12\right),
\left(\frac16,\frac56,\frac12\right).
\end{equation}
In each case,
since $\ul\a^{(1)}$ and $\ul\a$ agree up to permutation, 
the $F$-isocrystal $V^\crys(\ul \a,\sigma)$
is of rank $3$,
\[
V_{\ul\a}=\bigoplus_{i=0}^2K[t,(t-t^2)^{-1}]D^i\omega_{\ul\a}
\]

The following is the main result in this section.
\begin{thm}\label{m.thm-1}
Suppose $p>3$.
Let $a\in \Z_p$ satisfy $a(1-a)\not\equiv0$ mod $p$.
Put $b=\frac12(1-\sqrt{1-a})$.
Let $E_{\ul\a,s}$ be the elliptic curves as in \eqref{Ea-defn}.
Let $1-a_{p^2}(E_{\ul\a,b})T+p^2T^2$ be the characteristic polynomial of 
the $p^2$-th Frobenius
$\phi_{p^2}^{-1}\in G_{\F_{p^2}}$, namely $a_{p^2}(E_{\ul\a,b})\in \Z$ satisfies
\[
1-a_{p^2}(E_{\ul\a,b})+p^2=\sharp E_{\ul\a,b}(\F_{p^2}).
\]
Let
\[
d_{\ul\a}=\begin{cases}
-1&\ul\a=(\frac12,\frac12,\frac12),(\frac16,\frac56,\frac12)\\
-2&\ul\a=(\frac14,\frac34,\frac12)\\
-3&\ul\a=(\frac13,\frac23,\frac12)
\end{cases}
\]
and put
\[
A_{a,p}=\begin{cases}
a_{p^2}(E_{\ul\a,b})&\sqrt{1-a}\in \Z_p\\
\left(\frac{d_{\ul\a}}{p}\right)a_{p^2}(E_{\ul\a,b})
&\sqrt{1-a}\not\in \Z_p,\,E_{\ul\a,b}\text{: ordinary at $p$}\\
2p&\sqrt{1-a}\not\in \Z_p,\,E_{\ul\a,b}\text{: supersingular at $p$}.
\end{cases}
\]
Let $\sigma_a(t)=a^{1-p}t^p$.
Then
\[
\det(I-\Frob^{(\sigma_a)}_{\ul\a}|_{t=a}T)=\left(1-\left(\frac{1-a}{p}\right)pT\right)(1-A_{a,p}T+p^2T^2).
\]
\end{thm}
By isogeny
\eqref{isogeny1}, \eqref{isogeny2}, \eqref{isogeny3} and \eqref{isogeny4} in below,
one finds $a_{p^2}(E_{\ul\a,b})=a_{p^2}(E_{\ul\a,1-b})$, and hence $A_{a,p}$ does not depend on the choice of $b$.

\medskip

The proof of Theorem \ref{m.thm-1} shall be given in \S \ref{m.thm-1-pf-sect}.

\subsection{Frobenius on the elliptic curve $E_{\ul\a,s}$}\label{Frob-E-sect}
Let $s=\frac12(1-\sqrt{1-t})$, and $S:=\Spec W[s,t,(t-t^2)^{-1}]\to T$ a finite etale 
morphism
of degree $2$.
Let $\sigma$ be the $p$-th Frobenius given by $t^\sigma=ct^p$
as in \S \ref{F-H-sect}.
This extends on $\O(S)^\dag$ as follows,
\begin{align}
\sigma(\sqrt{1-t})&=\sqrt{1-t^\sigma}\\
&=\sqrt{1-t}\cdot(1-t)^{\frac{p-1}{2}}\sqrt{\frac{1-t^\sigma}{(1-t)^p}}\\
&=\sqrt{1-t}\cdot(1-t)^{\frac{p-1}{2}}\left(1+\frac12 pw(t)-\frac18(pw(t))^2+\cdots\right)
\end{align}
where $(1-t^\sigma)/(1-t)^p=1+pw(t)$.
For $\alpha\in W$ such that $(1-2\alpha)\alpha(1-\alpha)\not\equiv 0$ mod $p$,
one can define the evaluation at $s=\alpha$ 
of an element of $K[t,s,(t-t^2)^{-1}]^\dag$ or 
$\Q\ot(W[t,s,(t-t^2)^{-1}]^\wedge)$ the $p$-adic completion.
For example, if $t=a\in \Z_p$ satisfies $t^\sigma|_{t=a}=a$, then 
for $b=\frac12(1-\sqrt{1-a})$
\begin{equation}\label{1-b}
s^\sigma|_{s=b}=\begin{cases}
b&\sqrt{1-a}\in \Z_p\\
1-b&\sqrt{1-a}\not\in \Z_p.
\end{cases}
\end{equation}
Let
\[
g:\cE_{\ul\a}\lra S
\]
be the elliptic fibration whose general fiber is $E_{\ul\a,s}$.
We denote the fiber at $s=b$ by $E_{\ul\a,b}$, and
write $\cE_{\ul\a,K}=\cE_{\ul\a}\times_WK$, $\cE_{\ul\a,\ol\F_p}=\cE_{\ul\a}\times_W\ol\F_p$ etc.
Let
\[
H^\bullet_\rig(\cE_{\ul\a,\ol\F_p}/T_{\ol\F_p})
\]
be the rigid cohomology (cf. \cite{LS}).
The $\sigma$-linear $p$-th Frobenius endomorphism $\Phi_\cE$
is defined on the rigid cohomology, and hence on 
\[
K[t,s,(t-t^2)^{-1}]^\dag\ot_{\O(S_K)}
H^\bullet_\dR(\cE_{\ul\a,K}/S_K),
\]
thanks to the comparison theorem.
Since $\sigma(s)=4^{p-1}cs^p
+\cdots\in s^pW[[s]]$, the Frobenius
$\sigma$ extends on $W[[s]]=W[[t]]$ in a natural way, which we write by
the same notation.
The $p$-th Frobenius
$\Phi_\cE$ extends on $W[[s]]\ot_{\O(S)}H^1_\dR(\cE_{\ul\a,K}/S_K)$ as well. 
For each $\ul\a$,
we choose a Weierstrass equation
\begin{equation}\label{W-eq}
E_{\ul\a,s}:Y^2=4X^3-g_2(s)X-g_3(s)
\end{equation}
such that $g_2(s),g_3(s)\in W[s]$ satisfy $\Delta:=g_2^3-27g_3^2\in sW[s]$ and $g_2(0)g_3(0)\ne0$.
We then put
\[
\omega_{\cE,\ul\a}=\frac{dX}{Y},\quad
\eta_{\cE,\ul\a}=\frac{XdX}{Y}.
\]
\begin{lem}\label{ell-GM}
Let $D_s=s\frac{d}{ds}$ and put $E:=2g_2D_sg_3-3D_sg_2\cdot g_3$. Then
\begin{align*}
D_s\omega_{\cE,\ul\a}&=-\frac{D_s\Delta}{12\Delta}\omega_{\cE,\ul\a}
+\frac{3E}{2\Delta}
\eta_{\cE,\ul\a} ,\\
D_s\eta_{\cE,\ul\a}&=
-\frac{g_2E}{8\Delta}\omega_{\cE,\ul\a}
+\frac{D_s\Delta}{12\Delta}\eta_{\cE,\ul\a}.
\end{align*}
\end{lem}
\begin{pf}
This is well-known, e.g. \cite[Theorem 7.1]{A-C}.
\end{pf}
\begin{lem}\label{ell-GM-lem}
Let $\cD_S=K[s,t,\frac{d}{dt}]$ be the Wyle algebra of $S$. 
Let $D_s:=s\frac{d}{ds}$.
Then for $\ul\a=(\alpha_0,\alpha_1,\frac12)$, we have isomorphisms
\begin{equation}\label{ell-GM-lem-eq1-1}
\cD_S/\cD_SP_{\a_0\a_1}\os{\cong}{\lra}H^1_\dR(\cE_{\ul\a,K}/S_K),
\quad 1_{\cD_S}\longmapsto \omega_{\cE,\ul\a},
\end{equation}
\begin{equation}\label{ell-GM-lem-eq1}
\O(S_K)\ot_{\O(T)}V_{\ul \a}=\cD_S\cdot\omega_{\ul\a}\os{\cong}{\lra} 
\Sym^2H^1_\dR(\cE_{\ul\a,K}/S_K),
\quad \omega_{\ul\a}\longmapsto \omega_{\cE,\ul\a}^2
\end{equation}
of left $\cD_S$-modules, see \eqref{PHG} for the definition of $P_{\a_0\a_1}$
and \eqref{VHG} for $V_{\ul\a}$.
\end{lem}
\begin{pf}
We may replace $K$ with $\C$. It is a standard exercise to show that
there is a homology cycle $\gamma\in H_1(E_{\ul\a,s},\Q)$
that is a vanishing cycle at $s=0$, such that
\begin{equation}\label{ell-GM-lem-eq2}
\int_\gamma\omega_{\cE,\ul\a}=2\pi i\,{}_2F_1\left({\alpha_0,\alpha_1\atop 1};s\right).
\end{equation}
This implies $\cD_S\omega_{\cE,\ul\a}\subsetneq H^1_\dR(\cE_{\ul\a,K}/S_K)$, and hence
$\cD_S\omega_{\cE,\ul\a}=0$ as $H^1_\dR(\cE_{\ul\a,K}/S_K)$ is irreducible.
This shows that
the arrow \eqref{ell-GM-lem-eq1-1} is well-defined.
Since both side of \eqref{ell-GM-lem-eq1-1} are irreducible, 
it turns out to be bijective.

Next we show \eqref{ell-GM-lem-eq1}. Employing Clausen's formula \eqref{Clausen}
together with the transformation formula \cite[15.8.18]{NIST}, we have
\begin{equation}\label{Clausen-2}
{}_3F_2\left({\alpha_0,\alpha_1,\frac12\atop 1,1};t\right)={}_2F_1\left({\alpha_0,\alpha_1\atop 1};s\right)^2,\quad s=\frac12(1-\sqrt{1-t}).
\end{equation}
Therefore  we have
\[
\int_{\gamma\ot\gamma}P_{\ul\a}(\omega_{\cE,\ul\a}^2)=0
\]
and hence that the map
\[
\O(S)\ot_{\O(T)} V_{\ul\a}\cong\cD_S/\cD_SP_{\ul\a}\lra \Sym^2H^1_\dR(\cE_{\ul\a,K}/S_K),
\quad A\longmapsto A(\omega_{\cE,\ul\a}^2)
\]
is well-defined.
Since both sides are irreducible, this is bijective.
\end{pf}

Let us describe the Frobenius $\Phi_\cE$.
Let $W((t))^\wedge$ be the $p$-adic completion of $W((t))=W[[t]][t^{-1}]$.
Let $\ul\a=(\a_0,\a_1,\frac12)$.
According to \cite[\S 5.1]{AM}, there is a basis (de Rham symplectic basis)
\[
\wh\omega_{\ul\a,\cE},\,\wh\eta_{\ul\a,\cE}
\in 
W((t))^\wedge\ot_{\O(S)}
H^1_\dR(\cE_{\ul\a,K}/S_K)
\]
which satisfies (cf. \cite[Propositions 4.1, 4.2]{A-C})
\begin{equation}\label{wh-eq1}
\begin{pmatrix}
D_s(\wh\omega_{\cE,\ul\a})&
D_s(\wh\eta_{\cE,\ul\a})
\end{pmatrix}
=\begin{pmatrix}
\wh\omega_{\cE,\ul\a}
&\wh\eta_{\cE,\ul\a}
\end{pmatrix}
\begin{pmatrix}
0&0\\ \frac{D_sq}{q}&0\end{pmatrix}
\end{equation}
where $q\in W[[s]]$ is defined by $j(E_{\ul\a,s})=1/q+744+196884q+\cdots$, and
\begin{equation}\label{wh-eq2}
\begin{pmatrix}
\Phi_\cE(\wh\omega_{\cE,\ul\a})&
\Phi_\cE(\wh\eta_{\cE,\ul\a})
\end{pmatrix}
=\begin{pmatrix}
\wh\omega_{\cE,\ul\a}
&\wh\eta_{\cE,\ul\a}
\end{pmatrix}
\begin{pmatrix}
p&0\\ -p\tau^{(\sigma)}(s)&1\end{pmatrix}
\end{equation}
where $\tau^{(\sigma)}(s)\in W[[s]]$ is the power series given in \cite[Proposition 4.2]{A-C}.
The explicit relation with the basis $\{\omega_{\cE,\ul\a},\eta_{\cE,\ul\a}\}$ is
\begin{equation}\label{wh-eq3}
\wh\omega_{\cE,\ul\a}=\frac{\kappa}{F_{\a_0,\a_1}(s)}\omega_{\cE,\ul\a}
,\quad \wh\eta_{\cE,\ul\a}=-H(s)\omega_{\cE,\ul\a}
+\kappa^{-1}F_{\a_0,\a_1}(s)\eta_{\cE,\ul\a},\quad
\kappa:=\sqrt{-\frac{g_2(0)}{18g_3(0)}}
\end{equation}
where $H(s)$ is the power series given in \cite[Propositions 7.3]{A-C}.

\medskip

We pick up the results which shall be necessary in our proof of Theorem \ref{m.thm-1}.
\begin{prop}\label{Frob-ell}
Let
\begin{equation}\label{Frob-ell-eq1}
\begin{pmatrix}
\Phi_\cE(\omega_{\cE,\ul\a})&
\Phi_\cE(\eta_{\cE,\ul\a})
\end{pmatrix}
=\begin{pmatrix}
\omega_{\cE,\ul\a}
&\eta_{\cE,\ul\a}
\end{pmatrix}
\begin{pmatrix}
pA_0&B_0\\ pA_1&B_1\end{pmatrix}.
\end{equation}
Then all of $A_0,\ldots,B_1$ belong to $W[s,t,\sqrt{1-t}]^\dag$ and they satisfy
\begin{equation}\label{Frob-ell-eq2}
A_0B_1-A_1B_0=1
\end{equation}
and
\begin{equation}\label{Frob-ell-eq3}
B_1\equiv\epsilon\frac{F_{\a_0\a_1}(s)}{F_{\a_0\a_1}(s^\sigma)}\mod p,\quad \epsilon:=F(\kappa)/\kappa=\pm1.
\end{equation}
\end{prop}
\begin{pf}
See \cite[Propositions 4.3]{A-C}. Note that \eqref{Frob-ell-eq2} and \eqref{Frob-ell-eq3}
directly follow from the explicit descriptions of $A_0,\ldots,B_1$ in loc.cit.
\end{pf}
\subsection{Comparison Theorem}
\begin{thm}\label{ell-F-lem2}
Under the isomorphism
\begin{equation}\label{ell-F-lem2-eq1}
\O(S_K)^\dag\ot_{\O(T_K)}V_{\ul \a}
\cong \O(S_K)^\dag\ot_{\O(S_K)}\Sym^2H^1_\dR(\cE_{\ul\a,K}/S_K),
\end{equation}
induced from \eqref{ell-GM-lem-eq1}, the Frobenius $\Phi$ on the left
agrees with $\Sym^2\Phi_\cE$ on the right. 
\end{thm}
\begin{pf}
It follows from the main theorem of \cite{Dwork-unique} that they agree up to
a constant $c$, namely
$\Phi=c\,\Sym^2\Phi_\cE$.
We show $c=1$. To do this we employ the basis $\{\wh\omega_{\cE,\ul\a},\wh\eta_{\cE,\ul\a}\}$.
Since
\[
\langle\wh\xi_{\ul\a},\wh\eta_{\ul\a}\rangle=
\langle D\wh\omega_{\ul\a},D^2\wh\omega_{\ul\a}\rangle
\subset W((t))^\wedge\ot_{\O(T)}V_{\ul\a}
\]
by \eqref{FrobX-eq1},
one finds that
under the isomorphism \eqref{ell-F-lem2-eq1},
 this corresponds to
the subspace
\[
\langle D(\wh\omega^2_{\cE,\ul\a}),D^2(\wh\omega_{\cE,\ul\a}^2)\rangle
=\langle\wh\omega_{\cE,\ul\a}\wh\eta_{\cE,\ul\a},\wh\eta_{\cE,\ul\a}^2\rangle
\subset W((t))^\wedge\ot_{\O(S)}
\Sym^2H^1_\dR(\cE_{\ul\a,K}/S_K)
\]
where the above equality follows from \eqref{wh-eq1}.
Hence the isomorphism \eqref{ell-F-lem2-eq1} induces
an isomorphism
\[
\xymatrix{
(W((t))^\wedge\ot_{\O(T)}V_{\ul\a})/\langle\wh\xi_{\ul\a},\wh\eta_{\ul\a}\rangle
\ar[r]^-\cong& (W((t))^\wedge\ot_{\O(S)}
\Sym^2H^1_\dR(\cE_{\ul\a,K}/S_K))/
\langle\wh\omega_{\cE,\ul\a}
\wh\eta_{\cE,\ul\a},\wh\eta_{\cE,\ul\a}^2\rangle\\
W((t))^\wedge\,\wh\omega_{\ul\a}\ar@{=}[u]&
W((t))^\wedge\,\wh\omega_{\cE,\ul\a}^2\ar@{=}[u]
}
\]
compatible with respect to $\Phi$ and $c\,\Sym^2\Phi_\cE$.
It sends
$\wh\omega_{\ul\a}$ to $\kappa^2\wh\omega_{\cE,\ul\a}^2$ by \eqref{F-H-eq1} and \eqref{wh-eq3}
together with \eqref{Clausen-2}.
Since
\begin{align*}
\Phi(\wh\omega_{\ul\a})&\equiv
p^2\wh\omega_{\ul\a}\mod\langle\wh\xi_{\ul\a},\wh\eta_{\ul\a}\rangle
&\text{by \eqref{FrobX-eq6}}\\
\Sym^2\Phi_\cE(\wh\omega^2_{\cE,\ul\a})&\equiv p^2\wh\omega^2_{\cE,\ul\a}\mod
\langle\wh\omega_{\cE,\ul\a}\wh\eta_{\cE,\ul\a},\wh\eta^2_{\cE,\ul\a}\rangle&\text{by \eqref{wh-eq2}}
\end{align*}
one has
\[
c=\frac{F(\kappa^2)}{\kappa^2}=1
\]
as required.
\end{pf}
\begin{cor}\label{cor-3}
Let $\Frob_{\ul\a}^{(\sigma)}$ be the matrix defined in Definition \ref{Vcrys}, and
let $A_0,\ldots,B_1$ be as in \eqref{Frob-ell-eq1}.
Let $X$ be the $3\times3$ matrix defined by
\[
\begin{pmatrix}\omega_{\ul\a}&D\omega_{\ul\a}&D^2\omega_{\ul\a}
\end{pmatrix}
=\begin{pmatrix}\omega^2_{\cE,\ul\a}&
\omega_{\cE,\ul\a}\eta_{\cE,\ul\a}&
\eta^2_{\cE,\ul\a}
\end{pmatrix}
X
\]
under the identification \eqref{ell-F-lem2-eq1}.
Then 
\[
\Frob_{\ul\a}^{(\sigma)}=X^{-1}\begin{pmatrix}
p^2A_0^2&pA_0B_0&B_0^2\\
2p^2A_0A_1&p(A_0B_1+A_1B_0)&2B_0B_1\\
p^2A_1^2&pA_1B_1&B_1^2
\end{pmatrix}X^\sigma
\]
\end{cor}
\begin{pf}
Immediate from Theorem \ref{ell-F-lem2}.
\end{pf}

\subsection{Proof of Theorem \ref{m.thm-1}}\label{m.thm-1-pf-sect}
\begin{lem}\label{thm-4}
Let $q=p^m$, and
let $\phi_{q}=(\phi_p)^m\in G_{\F_p}$ be the $q$-th Frobenius.
Let $a\in W(\F_q)$ satisfy that $a(1-a)\not\equiv 0$ mod $p$.
Put
$b=\frac12(1-\sqrt{1-a})$.
Let $\sigma_a$ be given by $t^{\sigma_a}=F(a)a^{-p}t^p$ where $F$ is the $p$-th Frobenius on
$W(\F_q)$.
Then
\begin{equation}\label{thm-4-eq1}
\det(I-(\Frob_{\ul\a}^{(\sigma_a)})^{2m}|_{t=a}T)
=\det(1-(\phi^2_{q})^{-1}T\mid \Sym^2H^1_\et(E_{\ul\a,b,\ol\F_p},\Q_l))
\end{equation}
where $G_{\F_{q^2}}$ acts on $H^\bullet_\et(E_{\ul\a,b,\ol\F_p},\Q_l)$
as the reduction of $E_{\ul\a,b}$ is defined over $\F_{q^2}$.
If $\sqrt{1-a}\in W(\F_q)$, then
\begin{equation}\label{thm-4-eq2}
\det(I-(\Frob_{\ul\a}^{(\sigma_a)})^m|_{t=a}T)
=\det(1-(\phi_q)^{-1}T\mid \Sym^2H^1_\et(E_{\ul\a,b,\ol\F_p},\Q_l))
\end{equation}
where $G_{\F_q}$ acts on $H^\bullet_\et(E_{\ul\a,b,\ol\F_p},\Q_l)$
as the reduction of $E_{\ul\a,b}$ is defined over $\F_q$.
\end{lem}
\begin{pf}
Let $K_0$ be the fractional field of $W(\F_q)$, and we denote by $V^\dag_{\ul\a,K_0}$
the module \eqref{dagger-mod} replacing $K$ with $K_0$.
Then $\Phi|_{t=a}$ acts on $V^\dag_{\ul\a,K_0}|_{t=a}$ as a $F$-linear map, and 
the $m$-fold composition $\Phi^m|_{t=a}$ is a  $K_0$-linear map providing 
the characteristic polynomial of $(\Frob_{\ul\a}^{(\sigma_a)})^m$.
On the other hand, $\Phi^2_\cE|_{s=b}$ acts on $H^1_\dR(E_{\ul\a,b}/K_0)$ as a $F$-linear map, 
and the $K_0$-linear map $\Phi^{2m}_\cE|_{s=b}$ provides the right hand side of \eqref{thm-4-eq1}.
If $\sqrt{1-a}\in W(\F_q)$, then $\Phi_\cE|_{s=b}$ acts on $H^1_\dR(E_{\ul\a,b}/K_0)$, and
$\Phi^m_\cE|_{s=b}$ provides the right hand side of \eqref{thm-4-eq2}.
Therefore the lemma follows by Theorem \ref{ell-F-lem2}.
\end{pf}
\begin{lem}
Let $a\in \Z_p$ satisfy that $a(1-a)\not\equiv 0$ mod $p$, and $b:=\frac12(1-\sqrt{1-a})$.
Then
\begin{equation}\label{thm-2-eq1}
\det\Frob_{\ul\a}^{(\sigma)}|_{s=b}=\left(\frac{1-a}{p}\right)p^3.
\end{equation}
Let $A_0,\ldots,B_1$ be as in \eqref{Frob-ell-eq1}. 
If $\sqrt{1-a}\not\in \Z_p$, then
\begin{equation}\label{thm-2-eq2}
\Tr\,\Frob_{\ul\a}^{(\sigma_a)}|_{t=a}=\begin{cases}
(pA_0-B_1)^2|_{s=b}+p&\ul\a=(\frac12,\frac12,\frac12),(\frac16,\frac56,\frac12)\\
(p(A_0+2A_1)-B_1)^2|_{s=b}+p&\ul\a=(\frac13,\frac23,\frac12)\\
(p(A_0-\frac13A_1)-B_1)^2|_{s=b}+p&\ul\a=(\frac14,\frac34,\frac12).
\end{cases}
\end{equation}
\end{lem}
\begin{pf}
One has $\det\Frob_{\ul\a}^{(\sigma)}=p^3(1-2s)^3/(1-2s^\sigma)^3$ by
\eqref{Frob-ell-eq2} and  Corollary \ref{cor-3}. So \eqref{thm-2-eq1} is immediate by \eqref{1-b}.
To show \eqref{thm-2-eq2},
we need an explicit computation
of the matrix $X$.
For example, in case $\ul\a=(\frac12,\frac12,\frac12)$, one has
\[
X=\begin{pmatrix}1&-\frac13&\frac{7 s^{2}-7 s+1}{18(2s-1)^2}\\
0&\frac{1}{2s-1}&\frac{-10 s^{2}+10 s-1}{3(2s-1)^3}\\
0&0&\frac{1}{2(2s-1)^2}\end{pmatrix}
\]
by Lemma \ref{ell-GM}, and then \eqref{thm-2-eq2} follows. 
For other $\ul\a$'s, the proof is similar (left to the reader).
\end{pf}

\noindent{\it Proof of Theorem \ref{m.thm-1}}.
We want to compute the characteristic polynomial
\[
\det(I-\Frob^{(\sigma_a)}_{\ul\a}|_{t=a} T)=(1-\gamma_1T)(1-\gamma_2T)
(1-\gamma_3T)\in \Z[T].
\]
If $\sqrt{1-a}\in\Z_p$, then Theorem \ref{m.thm-1} is immediate from \eqref{thm-4-eq2}.
We suppose $\sqrt{1-a}\not\in\Z_p$ until the end of the proof.
We have 
\begin{equation}\label{proof-eq2}
\gamma_1\gamma_2\gamma_3
=-p^3
\end{equation}
by \eqref{thm-2-eq1}.
Let 
\[
\det(1-\phi_{p^2}^{-1}T\mid H^1_\et(E_{\ul\a,b,\ol\F_p},\Q_l))
=(1-\a_{p^2}(E_{\ul\a,b})T)(1-\beta_{p^2}(E_{\ul\a,b})T)\in\Z[T]
\]
be the characteristic polynomial
for the $p^2$-th Frobenius $\phi_{p^2}\in G_{\F_{p^2}}$.
We know from \eqref{thm-4-eq1} that the triplet $(\gamma^2_1,\gamma^2_2,\gamma^2_3)$
agrees with $(p^2,\alpha_{p^2}(E_{\ul\a,b})^2,\beta_{p^2}(E_{\ul\a,b})^2)$
up to permutation.
We may let  
\begin{equation}\label{proof-eq5}
\gamma_1=\pm p,\quad \gamma_2=\ve\a_{p^2}(E_b),\quad
\gamma_3=\ve'\beta_{p^2}(E_{\ul\a,b})
\end{equation}
where $\ve,\ve'=\pm1$.

If $E_{\ul\a,b}$ has a supersingular reduction at $p$, then $\ord_p(\alpha_{p^2}(E_{\ul\a,b}))>0$
and $\ord_p(\beta_{p^2}(E_b))>0$, and hence $\ord_p(\gamma_i)>0$ for every $i$.
By \eqref{thm-2-eq2}, there is $e\in W$ such that
\[
\gamma_1+\gamma_2+\gamma_3=e^2+p
\]
which is zero modulo $p$, whence
\begin{equation}\label{proof-eq3}
\gamma_1+\gamma_2+\gamma_3\equiv p\mod p^2.
\end{equation}
Suppose $\gamma_1=p$. Then $\gamma_2+\gamma_3\equiv 0$ mod $p^2$
by \eqref{proof-eq3}
and $\gamma_2\gamma_3=-p^2$ by \eqref{proof-eq2}. Since $\gamma_2+\gamma_3$ is a rational integer satisfying
$|\gamma_2+\gamma_3|\leq 2p$, it turns out that $\gamma_2+\gamma_3=0$ 
and hence $(\gamma_2,\gamma_3)=(p,-p)$ or $(-p,p)$. 
Suppose $\gamma_1=-p$. Then 
$\gamma_2+\gamma_3\equiv 2p$ mod $p^2$
and $\gamma_2\gamma_3=p^2$.
Since $|2p+np^2|\leq 2p$ $\Leftrightarrow$ $n=0$ as $p\geq 5$, 
one concludes $\gamma_2+\gamma_3=2p$ and hence $\gamma_2=\gamma_3=p$.
In both cases the triplet $(\gamma_1,\gamma_2,\gamma_3)$ agrees with
$(p,p,-p)$ up to permutation. This completes the proof of Theorem \ref{m.thm-1} when 
$E_{\ul\a,b}$ has a supersingular reduction.

The rest is the case that $E_b$ has an ordinary reduction.
We may let  $\ord_p(\alpha_{p^2}(E_{\ul\a,b}))=0$ and $\ord_p(\beta_{p^2}(E_{\ul\a,b}))>0$.
We first see that $\ve=\ve'$ in \eqref{proof-eq5}. Indeed, if $\ve\ne\ve'$, then
$\alpha_{p^2}(E_{\ul\a,b})-\beta_{p^2}(E_{\ul\a,b})$ is a rational integer
as $\gamma_1+\gamma_2+\gamma_3\in\Z$.
Since $\alpha_{p^2}(E_{\ul\a,b})+\beta_{p^2}(E_{\ul\a,b})\in \Z$, it turns out that
$\alpha_{p^2}(E_{\ul\a,b})$ is a rational integer satisfying $|\alpha_{p^2}(E_b)|=p$, and hence
$\alpha_{p^2}(E_{\ul\a,b})=\pm p$. This contradicts with the assumption 
$\ord_p(\alpha_{p^2}(E_{\ul\a,b}))=0$.
We thus have $\ve=\ve'$.
We want to show 
\begin{equation}\label{proof-eq4}
\ve=\left(\frac{d_{\ul\a}}{p}\right),
\end{equation}
which finishes the proof of Theorem \ref{m.thm-1} by \eqref{proof-eq2}.
We denote $f(s)_{<n}=\sum_{i<n}a_is^i$ the truncated polynomial
for a series $f(s)=\sum_{i=0}^\infty a_is^i$.
It follows from \eqref{Frob-ell-eq3} and \eqref{thm-2-eq2} that we have
\[
\ve\alpha_{p^2}(E_{\ul\a,b})\equiv 
\left(\frac{F_{\a_0\a_1}(s)}{F_{\a_0\a_1}(s^{\sigma})}\right)^2
\bigg|_{s=b}
\equiv (F_{\a_0\a_1}(s)_{<p})^2|_{s=b}
\mod p
\]
where
the second congruence follows from the Dwork congruence
(\cite[Theorem 3]{Dwork-p-cycle}).
On the other hand,
Dwork's unit root formula (\cite[(7.14)]{Put})
yields
\begin{align*}
\alpha_{p^2}(E_{\ul\a,b})
&=\frac{F_{\a_0\a_1}(s)}{F_{\a_0\a_1}(s^{\sigma^2})}
\bigg|_{s=b}
=
\frac{F_{\a_0\a_1}(s)}{F_{\a_0\a_1}(s^{\sigma})}
\frac{F_{\a_0\a_1}(s^\sigma)}{F_{\a_0\a_1}(s^{\sigma^2})}
\bigg|_{s=b}\\
&\os{\eqref{1-b}}{=}
\frac{F_{\a_0\a_1}(s)}{F_{\a_0\a_1}(s^{\sigma})}\bigg|_{s=b}
\cdot\frac{F_{\a_0\a_1}(s)}{F_{\a_0\a_1}(s^{\sigma})}
\bigg|_{s=1-b}\\
&\equiv F_{\a_0\a_1}(s)_{<p}|_{s=b}\cdot F_{\a_0\a_1}(s)_{<p}|_{s=1-b}\mod p
\end{align*}
Therefore we have
\[
\ve
\equiv F_{\a_0\a_1}(s)_{<p}|_{s=b}\cdot (F_{\a_0\a_1}(s)_{<p}|_{s=1-b})^{-1}
\mod p.
\]
We compute the right hand side.
To do this, we see the
Cartier operator $C$ for the reductions of 
$E_{\ul\a,b}$ and
$E_{\ul\a,1-b}$ respectively.
Let $\ul\a=(\frac12,\frac12,\frac12)$.
Recall that the elliptic curve
$E_{\ul\a,s}$ is defined by the Weierstrass equation $y^2=x(x-1)(x-s)$.
Then one has
\[
C\left(\frac{dx}{y}\right)=(-1)^{\frac{p-1}2}(F_{\frac12,\frac12}(s)_{<p}|_{s=b})\frac{dx}{y},\quad
C\left(\frac{dx}{y}\right)=(-1)^{\frac{p-1}2}(F_{\frac12,\frac12}(s)_{<p}|_{s=1-b})\frac{dx}{y}
\]
for $E_{\ul\a,b}$ and
$E_{\ul\a,1-b}$
respectively.
Using an isomorphism
\begin{equation}\label{isogeny1}
E_{\frac12,\frac12,\frac12,b}:y^2=x(x-1)(x-b)\lra 
E'_{\frac12,\frac12,\frac12,1-b}:-y^2=x(x-1)(x-1+b) 
\end{equation}
\[
(x,y)\longmapsto(1-x,y),
\]
it follows
\[
\ve\equiv
\frac{F_{\frac12,\frac12}(s)_{<p}|_{s=b}}{F_{\frac12,\frac12}(s)_{<p}|_{s=1-b}}\equiv
(\sqrt{-1})^p/\sqrt{-1}
=\left(\frac{-1}{p}\right)\mod p.
\]
This completes the proof of \eqref{proof-eq4} in case $\ul\a=(\frac12,\frac12,\frac12)$.
The proof in the other cases goes in the same way with use of the isogeny
\begin{equation}\label{isogeny2}
E_{\frac13,\frac23,\frac12,b}:
y^2=y^2=x^3+(3x+4b)^2\lra 
E'_{\frac13,\frac23,\frac12,1-b}:-3y^2=y^2=x^3+(3x+4-4b)^2
\end{equation}
\[
(x,y)\longmapsto\left(\frac{x^3+12x^2+48bx+64b^2}{-3x^2},\frac{x^3-48bx-128b^2}{9x^3}y\right),
\]
\begin{equation}\label{isogeny3}
E_{\frac14,\frac34,\frac12,b}:y^2=x(x^2-2x+1-b)\lra 
E'_{\frac14,\frac34,\frac12,1-b}:-2y^2=x(x^2-2x+b) 
\end{equation}
\[
(x,y)\longmapsto\left(-\frac12x+1+\frac{b-1}{2x},\frac{y(x^2+b-1)}{4x^2}\right)
\]
\begin{equation}\label{isogeny4}
E_{\frac16,\frac56,\frac12,b}:y^2=4x^3-3x+1-2b
\lra E'_{\frac16,\frac56,\frac12,1-b}: -y^2=4x^3-3x-1+2b
\end{equation}
\[
(x,y)\longmapsto(-x,y).
\]
This completes the proof of Theorem \ref{m.thm-1}.

\subsection{Cases $\ul\a=(\frac12,\frac12,\frac12),
(\frac16,\frac56,\frac12)$}
\begin{thm}\label{m.thm-2}
Let $\ul\a$ be either of $(\frac12,\frac12,\frac12),
(\frac16,\frac56,\frac12)$.
Let $C_{\ul\a,t}$ be the elliptic curves in \eqref{Ca-defn}.
Then 
\[
\det(I-\Frob^{(\sigma_a)}_{\ul\a}|_{t=a})
=\left(1-\left(\frac{1-a}{p}\right)pT\right)
\left(1-\left(\frac{1-a}{p}\right)a_{p^2}(C_{\ul\a,a})T+p^2T^2\right)
\]
for any $a\in \Z_p$ such that $a(1-a)\not\equiv 0$ mod $p$.
\end{thm}
\begin{pf}
The proof goes in a similar way to that of Theorem \ref{m.thm-1} and it is simpler.
We sketch the outline. The details are left to the reader.

The key formula on the hypergeometric function is
\begin{align*}
{}_3F_2\left({\alpha_0,\alpha_1,\frac12\atop 1,1};t\right)
&={}_2F_1\left({\frac12\alpha_0,\frac12\alpha_1\atop 1};t\right)^2&\text{(Clausen)}\\
&=(1-t)^{-\a_0}{}_2F_1\left({\frac12\alpha_0,1-\frac12\alpha_1\atop 1};
\frac{t}{t-1}\right)^2&\text{(\cite[15.8.1]{NIST}).}
\end{align*}
Let $\ul\a=(\frac12,\frac12,\frac12)$.
Then
\[
{}_3F_2\left({\frac12,\frac12,\frac12\atop 1,1};t\right)=(1-t)^{-\frac12}
{}_2F_1\left({\frac14,\frac34\atop 1};\frac{t}{t-1}\right)^2.
\]
Since the ``${}_2F_1$" in the right is the period of the elliptic curve 
$C_{\ul\a,t}$,
one can show an isomorphism
\[
V_{\ul\a}\cong\Sym^2H^1_\dR(\cC_{\ul\a,t}/T)
\ot\O(T)\sqrt{1-t}
\]
of $F$-isocrystals
in a similar way to the proof of Lemma \ref{ell-GM-lem},
where $\O(T)\sqrt{1-t}$ denotes the minus part $H^0_\dR(S/T)^-$ with respect to
the involution in $\mathrm{Aut}(S/T)\cong\Z/2\Z$.
The rest is automatic.

Let $\ul\a=(\frac16,\frac56,\frac12)$.
Then
\[
{}_3F_2\left({\frac16,\frac56,\frac12\atop 1,1};t\right)=
{}_2F_1\left({\frac{1}{12},\frac{5}{12}\atop 1};t\right)^2
\]
and the period of $C_{\ul\a,t}$ is
\[
(1-t)^{-\frac14}
{}_2F_1\left({\frac{1}{12},\frac{5}{12}\atop 1};t\right).
\]
Therefore one can show
\[
V_{\ul\a}\cong\Sym^2H^1_\dR(\cC_{\ul\a,t}/T)
\ot\O(T)\sqrt{1-t}
\]
as well. 

\end{pf}

\section{Galois representations of K3 surfaces of Hypergeometric type}\label{K3-sect}
\subsection{Lemmas on Galois representations of Kummer surfaces}
For an abelian surface $A$ over a field $F$ of characteristic $0$,
the minimal resolution of the quotient $A/\langle-1\rangle$
is a K3 surface, which we write by $\mathrm{Km}(A)$ and call the Kummer K3 surface.
Then one finds that there is a natural bijection
\[
H^2_{\et,\tr}(\ol{\mathrm{Km}(A)},\Q_l)\os{\cong}{\lra}
H^2_{\et,\tr}(\ol{A},\Q_l)
\]
on the transcendental part $H^2_{\et,\tr}(\ol X,\Q_l):=H^2_{\et}(\ol X,\Q_l)/\NS(\ol X)\ot\Q_l$ of the $l$-adic cohomology.
\begin{lem}\label{G-lem}
Suppose that $F$ is a number field. 
\begin{itemize}
\item[$(1)$]
The $G_F$-representation $H^2_{\et}(\ol{\mathrm{Km}(A)},\Q_l)$ is semisimple.
\item[$(2)$]
Let $V_{\Q_l}\subset H^2_{\et}(\ol{\mathrm{Km}(A)},\Q_l)$ and 
$V'_{\Q_l}\subset H^2_{\et}(\ol{\mathrm{Km}(A')},\Q_l)$ be sub $G_F$-representations.
Suppose that there is a set $S$ of primes of $F$ with density $1$ such that
\begin{equation}\label{G-lem-eq1}
\mathrm{Tr}(\phi_{\wp}\mid V_{\Q_l})
=\mathrm{Tr}(\phi_{\wp}\mid V'_{\Q_l})
\end{equation}
for every $\wp\in S$ where $\phi_\wp$
is the Frobenius at $\wp$.
Then there is an isomorphism
$V_{\Q_l}\cong V'_{\Q_l}$
of $G_F$-representations.
\end{itemize}
\end{lem}
\begin{pf}
(1) follows from a theorem of Faltings which asserts that $H^1_\et(\ol A,\Q_l)$ 
is semisimple.
(2) follows from the fact that $V_{\Q_l}$ and $V'_{\Q_l}$ are semisimple (cf.
\cite[I, 2.3]{serre}).
\end{pf}

\begin{lem}\label{sing-lem}
Let $X$ be a K3 surface over a number field $F$ such that $\rho(\ol X)\geq 19$.
Then $\ol X$ is singular $($i.e. $\rho(\ol X)=20)$
if and only if 
$\dim_{\Q_l}H^2_\et(\ol X,\Q_l(1))^{G_{F'}}=20$
for some finite extension $F'/F$.
\end{lem}
If we admit the Tate conjecture for cycles on $X\times X$ of codimension $2$,
then we can drop the condition ``$\rho(\ol X)\geq 19$".
\begin{pf}
By the assumption, $\ol X$ has the Shioda-Inose structure, namely
there is an elliptic curve $E$ over a number field such that the Kummer surface 
$\mathrm{Km}(\ol E\times\ol E)$ is isogenous to $\ol X$.
Then $\ol X$ is singular if and only if $\ol E$ has a CM.
The isogeny induces a surjective map
\[
\Sym^2H^1_\et(\ol E,\Q_l)\lra H^2_\et(\ol X,\Q_l)/\NS(\ol X)\ot\Q_l
\] 
of $G_{L}$-representations for some finite extension $L/F$.
Suppose that $\dim_{\Q_l}H^2_\et(\ol X,\Q_l(1))^{G_{F'}}=20$ for some $F'/F$. 
We may assume $L\subset F'$.
We want to show $\rho(\ol X)=20$. Suppose that 
$\rho(\ol X)=19$. Then since the above map is bijective,
the fixed part $(\Sym^2H^1_\et(\ol E,\Q_l)\ot\Q_l(1))^{G_{F'}}$ is $1$-dimensional.
This is equivalent to that $\dim_{\Q_l}\End(H^1_\et(\ol E,\Q_l))^{G_{F'}}= 2$.
Since $\End(\ol E)\ot\Q_l\cong \End(H^1_\et(\ol E,\Q_l))^{G_{F'}}$ 
by a theorem of Faltings, 
this implies that
$\ol E$ has a CM and hence that $\ol X$
is singular. This is a contradiction. This completes the proof of the ``if" part.
The ``only if" part is obvious.
\end{pf}

\subsection{K3 surfaces of Hypergeometric type}\label{K3HG-sect}
Let
\[
f:\cX\lra T=\Spec \Q[t,(t-t^2)^{-1}]
\]
be a projective smooth family of K3 surfaces.
Let $\cD=\Q[t,(t-t^2)^{-1},\frac{d}{dt}]$ be the Wyle algebra of $T$.
Put $\ol X_t:=\cX\times_T\ol{\Q(t)}$ and set
\[
V_\dR(\cX/T)=\mathrm{Coim}[H^2_\dR(\cX/T)\to H^2_\dR(\ol X_t)/\NS(\ol X_t)\ot\ol{\Q(t)}]
\]
 a $\cD$-module
which is free of finite rank over $\O(T)$.
\begin{defn}\label{K3.HG}
We call $f$ {\rm of hypergeometric type $F_{\ul\a}(t)$} if the following conditions hold.
\begin{itemize}
\item[{\rm(i)}] 
$\rho(\ol X_t)=19$ or equivalently $V_\dR(\cX/T)$ is of rank $3$ over $\O(T)$.
\item[{\rm (ii)}]
There is an isomorphism
$V_\dR(\cX/T)\cong \cD/\cD P_{\ul\a}$ of left $\cD$-modules, where
$P_{\ul\a}$ is the hypergeometric differential operator \eqref{PHG}.
\end{itemize}
\end{defn}
Let $\ul\a$ be either of the following as before,
\[
\left(\frac12,\frac12,\frac12\right),
\left(\frac13,\frac23,\frac12\right),
\left(\frac14,\frac34,\frac12\right),
\left(\frac16,\frac56,\frac12\right).
\]
Let $p>3$ be a prime and put $W=W(\ol\F_p)$ and $K=\Frac W$.
Suppose that there is an integral regular flat model
\[
f_{\Z_{(p)}}:\cX_{\Z_{(p)}}\lra T_{\Z_{(p)}}
\]
over the ring $\Z_{(p)}\subset \Q$ such that $f_{\Z_{(p)}}$ is smooth projective. 
For $a\in \Z_{(p)}$ such that $a(1-a)\not\equiv0$ mod $p$, we denote by $X_a$ the fiber
at $t=a$.
We denote $\cX_W=\cX_{\Z_{(p)}}\times_{\Z_{(p)}}W$,
$\cX_K=\cX_W\times_WK$, \ldots as before.
Recall from Definition \ref{Vcrys} the $F$-isocrystal 
\[
V^\crys(\ul \a,\sigma)=\left(V_{\ul\a},
V_{\ul\a}^\dag,
\nabla,\Phi,\sigma\right).
\]
Thanks to the comparison
\[
K[t,(t-t^2)^{-1}]^\dag\ot_{\O(T_K)} H^2_\dR(\cX/T)
\cong
H^2_\rig(\cX_{\ol\F_p}/T_{\ol\F_p})
\]
of the de Rham and rigid cohomology, the $p$-th Frobenius on
\[
K[t,(t-t^2)^{-1}]^\dag\ot_{\O(T_K)} 
V_\dR(\cX_K/T_K)
\]
is defined, which we denote by $\Phi_\cX$.
We fix an isomorphism $V_\dR(\cX/T)\cong \cD/\cD P_{\ul\a}$ (this is unique
up to scalar as both sides are irreducible).
It follows from \cite{Dwork-unique} that there is a constant $\ve_p$ such that
\[
\Phi=\ve_p\Phi_\cX.
\]
Let $a\in \Z_{(p)}$ such that $a(1-a)\not\equiv0$ mod $p$ and 
let $\sigma=\sigma_a$
given by $\sigma_a(t)=a^{1-p}t^p$.
By Theorem \ref{m.thm-1}, the characteristic polynomial of
$\Phi_\cX|_{t=a}$ on $V_\dR(X_a/\Q_p)$ is
\[
\left(1-\left(\frac{1-a}{p}\right)p\ve_pT\right)(1-\ve_pA_{a,p}T+(\ve_ppT)^2).
\]
Since this is a polynomial with coefficients in $\Z$ and $p\pm A_{a,p}\ne0$, it turns out that $\ve_p\in \Q^\times$ (not depend on $a$).
Moreover since $(\det V_\dR(X_a/\Q_p))^{\ot2}$ is 
isomorphic to the Tate object $\Q_p(-6)$, one has $\ve_p^6=1$
and hence $\ve_p=\pm1$.
Using the $l$-adic cohomology, 
\[
V_\et(\ol X_a)_{\Q_l}:=\mathrm{Coim}[H^2_\et(\ol X_a,\Q_l)
\to H^2_\et(\ol X_t,\Q_l)/\NS(\ol X_t)\ot\Q_l],\quad
\ol X_a:=X_a\times_{\Z_{(p)}}\ol\Q,
\]
one has
\[
\det(\phi_p^{-1}\mid(\det V_\et(\ol X_a)_{\Q_l})\ot\Q_l(3))
=\ve^3_p\left(\frac{1-a}{p}\right)
=\ve_p\left(\frac{1-a}{p}\right).
\]
Define a quadratic character $\chi_{\cX/T}$ by
\begin{equation}\label{chiX}
\chi_{\cX/T}:G_\Q\lra \mathrm{Aut}((\det V_\et(\ol X_a)_{\Q_l})\ot\Q_l(3)\ot\chi_{1-a})
\cong\Q_l^\times
\end{equation}
where $\chi_{1-a}$ denotes the quadratic character for $\Q(\sqrt{1-a})$.
This does not depend on either $a$ or $l$, takes value in $\{\pm1\}$
and satisfies $\chi_{\cX/T}(\phi_p)=\ve_p$ for almost all $p$.

\medskip

Summing up the above, we have 
\begin{thm}\label{thm-1}
Let the notation and assumption be as above. Then
\[
\det(1-\phi_p^{-1}T\mid V_\et(\ol X_a)_{\Q_l})
=\left(1-\left(\frac{1-a}{p}\right)
\chi_{\cX/T}(\phi_p)pT\right)(1-\chi_{\cX/T}(\phi_p)A_{a,p}T+p^2T^2).
\]
\end{thm}

\begin{thm}\label{thm-2}
Let $F$ be a number field.
Let $a\in F\setminus\{0,1\}$ be arbitrary.
Then
there is an isomorphism
\[V_\et(\ol X_a)_{\Q_l}\cong \Sym^2H^1_\et(\ol E_{\ul\a,b},\Q_l)\ot\chi_{\cX/T}
\]
of $G_{F(\sqrt{1-a})}$-representations
where $E_{\ul\a,b}$ is the elliptic curve in Theorem \ref{m.thm-1}.
\end{thm}
\begin{pf}
By virtue of Lemma \ref{G-lem} (2), it is enough to show that  
\[
\Tr(\phi_\wp\mid V_\et(\ol X_a)_{\Q_l})
=\Tr(\phi_\wp\mid \Sym^2H^1_\et(\ol E_{\ul\a,b},\Q_l)\ot\chi_{\cX/T}).
\]
for almost all primes $\wp$ of $F(\sqrt{1-a})$.
However this follows from Lemma \ref{thm-4} \eqref{thm-4-eq2}.
\end{pf}

In the case 
$\ul\a=(\frac12,\frac12,\frac12)$ or $(\frac16,\frac56,\frac12)$,
the same discussion replacing Theorem \ref{m.thm-1} with Theorem \ref{m.thm-2} 
yields the following.
\begin{thm}\label{thm-3}
Let $\ul\a$ be either of $(\frac12,\frac12,\frac12)$ or $(\frac16,\frac56,\frac12)$.
Then, for $a\in F\setminus\{0,1\}$,
there is an isomorphism
\[V_\et(\ol X_a)_{\Q_l}\cong \Sym^2H^1_\et(\ol C_{\ul\a,a},\Q_l)\ot\chi_{\cX/T}\ot\chi_{1-a}
\]
of $G_F$-representations
where $C_{\ul\a,a}$ is the elliptic curve in Theorem \ref{m.thm-2}.
\end{thm}
\begin{thm}\label{s.thm-4}
Let $a\in \ol\Q\setminus\{0,1\}$. Then $X_a$ is singular (i.e. the Picard number $20$) if and only if
$E_{\ul\a,b}$ has a CM.
\end{thm}
\begin{pf}
Immediate from Theorem \ref{thm-2} and Lemma \ref{sing-lem}.
\end{pf}
It is known that there are 13 rational CM $j$-invariants, and
58 CM $j$-invariants that are quadratic numbers but not rational (\cite{D-LR}).
The complete list of them is available by SAGE.
We thus have the complete list of rational $a$'s such that $X_a$ is singular.
\bigskip
\begin{flushleft}
\begin{tabular}{|c|c|}\hline
&$a\in\Q\setminus\{0,1\}$ such that $X_a$ is singular\\
\hline
$\displaystyle\begin{matrix}\vspace{1cm}\end{matrix}$
$\displaystyle\ul\a=\left(\frac12,\frac12,\frac12\right)$&$\displaystyle-1,4,\frac14,-8,-\frac18,64,\frac{1}{64}$
\\
\hline
$\displaystyle\begin{matrix}\vspace{1cm}\end{matrix}$
$\displaystyle\ul\a=\left(\frac13,\frac23,\frac12\right)$&
 $\displaystyle
-4,\frac12,\frac{-1}{2^4},\frac{-1}{2^{10}},\frac{-3^2}{2^4},\frac{3^3}{2^4},
\frac{2}{3^3},\frac{3^3}{2},\frac{-1}{2^45},\frac{4}{5^3},\frac{-1}{2^45^6},\frac{-1}{2^33^37}
$\\
\hline$\displaystyle\begin{matrix}\vspace{1cm}\end{matrix}$
 $\displaystyle\ul\a=\left(\frac14,\frac34,\frac12\right)$&
 $\displaystyle
\frac{-1}{2^{2}},\,\frac{1}{3^{2}},\,
\frac{1}{3^{4}},\,
\frac{-1}{2^{4}3},\,
\frac{-2^4}{3^2},\,
\frac{-1}{2^{2}3^{4}},\,
\frac{2^5}{3^{4}},\,\frac{2^8}{3^4},\,
\frac{-1}{2^{6}3^{4}5},\,
\frac{1}{7^{4}},\,
\frac{-1}{2^{2}3^{4}7^{4}},\,
\frac{-2^{8}}{3^{4}7^{2}},\,
\frac{1}{3^{4}11^{2}},\,\frac{1}{3^{8}11^{4}}
$\\
\hline$\displaystyle\begin{matrix}\vspace{1cm}\end{matrix}$
$\displaystyle\ul\a=\left(\frac16,\frac56,\frac12\right)$&
 $\displaystyle
\frac{-1}{2^{9}},\,
\frac{-3^{3}}{2^{9}},\,
 \frac{2^{2}}{5^{3}},\,
\frac{-2^{6}}{5^{3}},\,\frac{-1}{2^{12}5^{3}},\,
\frac{-3^{2}}{2^{9}5^{3}},\,\frac{3^{3}}{5^{3}},\,
\frac{2^3}{11^3},\,
\frac{-1}{2^{9}5^{3}11^{3}},\,\frac{2^{6}}{5^{3}17^{3}},\,\frac{-1}{2^{12}5^{3}23^{3}29^{3}}$\\
      \hline
\end{tabular}
\end{flushleft}
\begin{rem}
It is natural to expect that the Shioda-Inose structure of 
$\ol X_a$ is given by $\mathrm{Km}(\ol E_{\ul\a,b}\times \ol E_{\ul\a,b})$,
the author does not have a proof though.
\end{rem}
\begin{rem}\label{chi-rem}
For the K3 families (i), \ldots, (iv) in \S \ref{intro-sect}, the character $\chi_{\cX/T}$
is trivial.
For example, let $\cX\to T$ be the family (iii).
There is an isomorphism (\cite[Corollary 3.6]{As-Ross1})
\[
\cD/\cD P_{\ul\a}\os{\cong}{\lra} V_\dR(\cX/T),\quad A\longmapsto
A\omega_{1,1,1}
\]
where $\omega_{i,j,k}$ is the regular $2$-form defined in \cite[(2.7)]{As-Ross1}.
Letting $\omega:=\omega_{1,1,1}$, we construct a basis
$\{\wh\omega,\wh\xi,\wh\eta\}$ in the same as in \S \ref{F-H-sect}.
Then it is shown that
\[
\Phi_\cX(\wh\omega)\equiv p^2\wh\omega\mod\langle\wh\xi,\wh\eta\rangle
\]
in the proof of \cite[Theorem 4.5]{As-Ross2}. This shows $\ve_p=1$, 
namely $\chi_{\cX/T}$ is trivial.
For the other K3 families, one can show that $\chi_{\cX/T}$ is trivial in the same way
(the argument in the proof of \cite[Theorem 4.5]{As-Ross2} works).
\end{rem}


\begin{thebibliography}{AAAA}


\bibitem[AOP]{AOP}
Ahlgren, S.; Ono, K.; Penniston, D.,
{\it Zeta Functions of an infinite family of K3 surfaces}.
American Journal of Mathematics, Vol. 124, No. 2 (Apr., 2002), pp. 353--368.


 \bibitem[As1]{As-Ross1}
Asakura, M.,
{\it A generalization of the Ross symbols in higher 
$K$-groups and hypergeometric functions I},
arXiv.2003.10652 


 \bibitem[As2]{As-Ross2}
\bysame,
{\it A generalization of the Ross symbols in higher 
$K$-groups and hypergeometric functions II},
arXiv:2102.07946.

 \bibitem[A-C]{A-C}
Asakura, M. and Chida, M.,
{\it A numerical approach toward the 
$p$-adic Beilinson conjecture for elliptic curves over $\Q$}.
arXiv:2003.08888
 \bibitem[A-M]{AM}
Asakura, M. and Miyatani, K.,
{\it Milnor $K$-theory, $F$-isocrystals and syntomic regulators}.
arXiv:2007.14255.


 \bibitem[D-LR]{D-LR}
Daniels, H-B., Lozano-Robledo, A.,
{\it On the number of isomorphism classes of CM elliptic curves defined over a number field.} 
J. Number Theory \textbf{157} (2015), 367--396.
\bibitem[Dw1]{Dwork-p-cycle}
Dwork, B.,
{\it $p$-adic cycles}.
Publ. Math. IHES, tome 37 (1969), 27--115.

\bibitem[Dw2]{Dwork-unique}
\bysame,
{\it On the uniqueness of Frobenius operator on differential equations}.
in Algebraic Number Theory - in Honor of K. Iwasawa, Advanced Studies in Pure Math. 17, Academic Press, Boston, 1989, 89--96.


\bibitem[E-S]{ES}
Elkies, N. and Sch\"utt, M.,
{\it K3 families of higher Picard rank}. preprint.

\bibitem[vG-T]{GT}
van Geemen, B.,  Top, J.,
{\it An isogeny of K3 surfaces.} Bull. London Math. Soc. \textbf{38} (2006), 
no. 2, 209--223. 
\bibitem[Go1]{Go1}
Goodson, H., 
{\it Hypergeometric functions and relations to Dwork hypersurfaces}. 
Int. J. Number Theory, Volume \textbf{3} (2), 2017, 439--485.
\bibitem[Go2]{Go2}
\bysame,
{\it A complete hypergeometric point count formula for Dwork hypersurfaces}. 
J. Number Theory, Volume \textbf{179}, 2017, 142--171.
\bibitem[Ka]{Katz}
Katz, N., {\it Another look at the Dwork family}, 
Algebra, arithmetic, and geometry: in honor of Yu. I. Manin. Vol. II, 89--126, Progr. Math., 
\textbf{270}, Birkhauser, 2009. 
\bibitem[Ko]{Koblitz}
Koblitz, N., {\it
The number of points on certain families of hypersurfaces over finite fields}. 
Compositio Math., tome {\bf 48}, 1983, 3--23.
\bibitem[Mc]{Mc}
McCarthy, D., 
{\it The number of $\F_q$-points on Dwork hypersurfaces and hypergeometric functions.} Res. Math. Sci., Volume \textbf{4} (4), 2017, 1--15.
\bibitem[Mi]{Miyatani}
K. Miyatani, K., 
{\it Monomial deformations of certain hypersurfaces and two hypergeometric functions.}
Int. J. Number Theory, Volume \textbf{11} (8), 2015, 2405--2430.
\bibitem[Mo]{morrison} Morrison, D. R.,
{\it On K3 surfaces with large Picard number}, 
Invent. Math. {\bf 75}, (1984), 105--121. 
\bibitem[Na]{nas}
Naskrecki, B.,
{\it On a certain hypergeometric motive of weight 2 and rank 3}.
arXiv:1702.07738
\bibitem[O]{otsubo} Otsubo, N.,
{\it Hypergeometric functions over finite fields}. arXiv:2108.06754. 
\bibitem[VdP]{Put}
Van der Put, M.,
{\it The cohomology of Monsky and Washnitzer.}
Introductions aux cohomologies $p$-adiques (Luminy, 1984). 
Mem. Soc. Math. France (N.S.)  No. 23  (1986), 33--59. 

\bibitem[Se]{serre}
Serre, J-P.,
{\it Abelian $l$-adic representations and elliptic curves.}
Benjamin, 1968.

\bibitem[LS]{LS}
Le Stum, B.,
{\it Rigid cohomology.} Cambridge Tracts in Mathematics, 172. Cambridge University Press, Cambridge, 2007. xvi+319 pp. 

\bibitem[NIST]{NIST} {\it NIST Handbook of Mathematical Functions. }
Edited by Frank W. J. Olver, Daniel W. Lozier, Ronald F. Boisvert and Charles W. Clark. 
Cambridge Univ. Press, 2010.
\end{thebibliography}
\end{document}